\documentclass[reqno, 12pt]{amsart}

\usepackage{amsmath,amssymb,amsthm,amsfonts}
\usepackage{fullpage}
\usepackage{xcolor}
\usepackage[colorlinks=true,linkcolor=blue,anchorcolor=blue,citecolor=red]{hyperref}
\allowdisplaybreaks % allow equations break between two pages
\usepackage{cite}
\usepackage{colonequals}

\newtheorem{theorem}{Theorem}[section]
\newtheorem{lemma}[theorem]{Lemma}
\newtheorem{corollary}[theorem]{Corollary}

\theoremstyle{definition}

\theoremstyle{remark}
\newtheorem{remark}[theorem]{Remark}
\newcommand{\Z}{\mathbb{Z}}

\newcommand{\R}{\mathbb{R}}

\renewcommand{\P}{\mathbb{P}}

\newcommand{\F}{\mathbb{F}}

\newcommand{\cA}{\mathcal{A}}
\newcommand{\cB}{\mathcal{B}}

\renewcommand{\mod}[1]{\,(\mathrm{mod}\,#1)}
\newcommand{\of}[1]{\left(#1\right)}
\newcommand{\set}[1]{\left\{#1\right\}}

\author{Huixi Li}
\address{School of Mathematical Sciences and LPMC, Nankai University, Tianjin 300071, China}
\email{lihuixi@nankai.edu.cn}

\author{Biao Wang}
\address{School of Mathematics and Statistics, Yunnan University, Kunming, Yunnan 650091, China}
\email{bwang@ynu.edu.cn}

\author{Chunlin Wang}
\address{School of Mathematical Sciences, Sichuan Normal University, Chengdu 610064, China}
\email{c-l.wang@outlook.com}

\author{Shaoyun Yi}
\address{School of Mathematical Sciences, Xiamen University, Xiamen, Fujian 361005, China}
\email{yishaoyun926@xmu.edu.cn}

\date{\today}

\makeatletter
\@namedef{subjclassname@2020}{\textup{}2020 Mathematics Subject Classification}
\makeatother

\title[]{On Erd\H{o}s covering systems in global function fields}

\begin{document}
	
\begin{abstract} 
A covering system of the integers is a finite collection of arithmetic progressions whose union is the set of integers. A well-known problem on covering systems is the minimum modulus problem posed by  Erd\H{o}s in 1950, who asked whether the minimum modulus in such systems with distinct moduli can be arbitrarily large. This problem was resolved by Hough in 2015, who showed that the minimum modulus is at most $10^{16}$. In 2022, Balister, Bollob\'as, Morris, Sahasrabudhe and Tiba reduced Hough's bound to $616,000$ by developing Hough's method. They call it the distortion method. In this paper, by applying this method, we mainly prove that there does not exist any covering system of multiplicity $s$ in any global function field of genus $g$ over $\mathbb{F}_q$ for $q\geq (1.14+0.16g)e^{6.5+0.97g}s^2$. In particular, there is no covering system of $\mathbb{F}_q[x]$ with distinct moduli for $q\geq 759$.
\end{abstract}

\maketitle

\section{Introduction and statement of results}\label{sec_intro}

A \textit{covering system} is a finite collection of arithmetic progressions $\{A_1,...,A_m\}$ that covers the integers, that is, $\bigcup_{i=1}^m A_i=\Z$. Covering systems were first introduced by Erd\H{o}s \cite{Erdos1950} over 70 years ago. And he posed many problems about these systems with distinct moduli.  One of the most famous problems, which is also Erd\H{o}s' favorite one \cite{Erdos1995}, is the so-called \textit{minimum modulus problem}: whether the minimum modulus in covering systems with distinct moduli can be arbitrarily large. In 2015, following the breakthrough by Filaseta, Ford, Konyagin, Pomerance and Yu \cite{FFSPY2007}, Hough \cite{Hough2015} finally resolved this problem and proved that the minimum modulus is at most $10^{16}$ in any covering system of the integers with distinct moduli. Building on Hough's ideas, Balister, Bollob\'as, Morris, Sahasrabudhe and Tiba \cite{BBMST2022} reduced Hough's bound to $616,000$ by developing a new method to estimate the density of the uncovered set. They called it the \textit{distortion method} in \cite{BBMST2020_2}. Recently, Cummings, Filaseta and Trifonov \cite{CFT2022} reduced this bound to 118 for  covering systems with distinct squarefree moduli. They also proved that the $j$-th smallest modulus in a minimal covering system with distinct moduli is bounded. Later, Klein, Koukoulopoulos and Lemieux \cite{KKL2022} gave a specified upper bound for the $j$-th smallest modulus. Following their work \cite{KKL2022}, the authors Li, Wang and Yi \cite{LiWangYi2023pre} developed the distortion method in the setting of finite sets and provided a solution to a generalization of Erd\H{o}s' minimum modulus problem in number fields. Moreover, the authors \cite{LiWangWangYi2023pre} proved an analogue of \cite[Theorem~3]{KKL2022} to polynomial rings over finite fields. 

In this paper, we continue to study the covering systems of polynomial rings over finite fields.   By Hough's result, we know that the minimum modulus of a covering system of the integers with distinct moduli is bounded. Let $\F_q$ be a finite field with $q$ elements. For any polynomial $f(x)\in \F_q[x]$, the norm of $f(x)$ is $|f|=q^{\deg f}$. An obvious observation is that if $q$ is sufficiently large, then $|f|$ would be sufficiently large, too. This implies that the smallest norm of a covering system of $\F_q[x]$ is sufficiently large if $q$ is sufficiently large. Thus, motivated by Hough's result, we believe that if a covering system of $\F_q[x]$ with distinct moduli exists, then $q$ is bounded. In this paper, we will mainly prove this statement. Indeed, we will prove it in the setting of global function fields. We first introduce some necessary notation on global function fields as follows.

Let $F=\mathbb{F}_q$. Let $K/F$ be a global function field over the field of constants $F$. Let $\mathcal{S}_K$ be the set of all primes $P$ of $K$, and let $S\subset \mathcal{S}_K$ be any non-empty, finite set of primes. Define
\[
O_S\colonequals\{a\in K\mid \mathrm{ord}_P(a)\geq 0, \forall P\notin S\},
\]
the ring of $S$-integers. Here, $\mathrm{ord}_P\colon K\to\Z$ is the order function attached to $P$; see \cite[Chapter~9]{AtiyahMacdonald1969}. By \cite[Theorem~14.5]{Rosen2002}, there exist elements $x\in K$ such that the poles of $x$ consist precisely of the elements of $S$. And for any such element $x$, the integral closure of $F[x]$ in $K$ is $O_S$. Moreover, $O_S$ is a Dedekind domain, and there is a one-to-one correspondence between the nonzero prime ideals of $O_S$ and the primes of $K$ not in $S$. More precisely, a prime in $K$ is the maximal ideal, also denoted by $P$, of the discrete valuation ring $O_P$ with $F\subset O_P$ and the quotient field of $O_P$ is just $K$. The norm $|P|$ of $P$ is the size of the residue field $\kappa_P$ of $O_P$, i.e., $|P|=[O_P\colon P]=\#\kappa_P$, which is a power $q^{\deg P}$ of the cardinality of the ground field $F$. Here the exponent $\deg P$ is called the degree of $P$. Furthermore, it follows from the well-known Riemann-Roch theorem (e.g., see \cite[Theorem~5.4]{Rosen2002}) that there is an integer $g\geq 0$ uniquely determined by $K$, which is called the
genus of $K$. See also \cite{DuanWangYi2021} for a nice introduction to various aspects of global function fields.

By the celebrated work of Weil \cite{Weil1948}, we have the Riemann Hypothesis for global function fields; see for example \cite[Theorem~5.10]{Rosen2002}. As a consequence, we have an analogue of the prime number theorem for global function fields with the best bound (e.g., see \cite[Theorem~5.12]{Rosen2002}) as follows:
\begin{equation}\label{eqn_PNT_for_FF}
\pi_K(n)\colonequals\#\{P\mid \deg P=n\}=\frac{q^n}{n}+O\left(\frac{q^{\frac{n}{2}}}{n}\right).
\end{equation}

Let $I=\prod_PP^{a_P}$ be an ideal of $O_S$ such that every $a_P=\mathrm{ord}_P(I)$ is a
non-negative integer, and $a_P=0$ for all but finitely many $P$. It is clear that the degree $\deg I$ of $I$ equals to $\sum_Pa_P\deg P$ and hence the norm $|I|$ of $I$ is $q^{\deg I}$. We say that a \textit{covering system} of $O_S$ is a finite collection  $\mathcal{A}=\{A_{i}\colon 1\leq i\leq m\}$ of congruences in $O_S$ such that $O_S=\bigcup_{i=1}^m(a_i+I_{i})$, where $A_i=a_i+I_{i}$ for some $a_i\in O_S$ and some ideal $I_{i}\subseteq O_S$. 
Let 
\[
    m(\cA)=\max_{ \text{Ideals } I \subset O_S}\#\{1\leq i\leq m\colon I_i=I\}
\]
be the multiplicity of $\cA$. Then the main result of our paper is stated as follows.

\begin{theorem}\label{mainthm}
Let $F=\mathbb{F}_q$ be a finite field. Let $K/F$ be a global function field over the field of constants $F$. Let $g$ be the genus of $K$. Let $S$ be a non-empty finite set of prime ideals, and let $O_S$ be the ring of $S$-integers.  Let $s\geq1$ be an integer. Then there does not exist any covering system of $O_S$ with multiplicity $s$ for $q\geq (1.14+0.16g)e^{6.5+0.97g}s^2$. 
\end{theorem}

In particular, when $O_S=\F_q[x]$ is a polynomial ring over $\F_q$, then we know that the genus $g=0$ in this case. If the moduli of a covering system is distinct, then the multiplicity $s=1$. As a corollary of Theorem~\ref{mainthm}, we arrive at the following result on the minimum modulus problem for covering systems over $\F_q[x]$. 

\begin{corollary}\label{maincor}
	There does not exist any covering system of $\mathbb{F}_q[x]$ of multiplicity $s$ for $q\geq 759s^2$. In particular, there is no covering system of $\mathbb{F}_q[x]$ with distinct moduli for $q\geq 759$.
\end{corollary}

\begin{remark}
    The constant $759$ is not optimal. One may numerically improve it a little bit further to 739 by rounding decimals in Lemmas~\ref{lem_mertens_estimates} and \ref{lem_2nd_moments_estimate} with precision 0.0001. In the proofs of our results, we round decimals with precision 0.01 for brevity. Furthermore, if replacing  the condition $q\geq 700$  by $q\geq724$ in Lemmas~\ref{lem_mertens_estimates} and \ref{lem_2nd_moments_estimate}, then one can obtain a better constant $725$ in Corollary~\ref{maincor}. This is the optimal numerical estimation in our method. Since $q$ is a prime power, we get that $q\neq 720, 721, \dots, 724$.  Thus, we conclude that  there is no covering system of $\mathbb{F}_q[x]$ with distinct moduli for $q\geq 720$.

    For $q=2$, the following are two examples of covering systems of $\mathbb{F}_2[x]$ with distinct moduli found by Azlin \cite{Azlin2011thesis} in his master's thesis: 
	\begin{align*}
		&(1). \quad  \{
0 \mod{x}, \quad
0 \mod{x + 1}, \quad
1 \mod{x^2 + x}
\};\\
&(2). \quad \{
0 \mod{x^2}, \quad
x \mod{x^2 + 1}, \quad
x + 1 \mod{x^2 + x}, \quad
1 \mod{x^3 + x}, \\
&\qquad\qquad x^2 + x \mod{x^3 + x^2}, \quad
x^3 + x^2 + x \mod{x^4 + x^2}
\}.
	\end{align*}
An interesting problem is to find examples of covering systems of $\mathbb{F}_q[x]$ with distinct moduli for $3\leq q <720$. As regards the integers $\Z$ case, while Nielsen \cite{Nielsen2009} constructed a covering system with minimum
modulus 40, Owens \cite{Owens2014} improved this record to 42.

\end{remark}

\begin{remark} Since the analogue \eqref{eqn_PNT_for_FF} of the prime number theorem for global function fields holds, one may generalize \cite[Theorem~1.4]{LiWangWangYi2023pre} as follows: let $\mathcal{A}=\{a_i+I_{i}\colon 1\leq i\leq m\}$ be covering system of $O_S$ of multiplicity $s$, then there exists some constant $c>0$ depending only on the genus $g$ of $K$ such that
\begin{equation}\label{eqn_min_degree}
	\min_{1\leq i\leq m}\deg I_i \leq 3(c+3\log_qs)\log(cs+3s\log_qs).
\end{equation}
This result gives a solution to a generalization of Erd\H{o}s' minimum modulus problem for global function fields in the degree aspect. It is also an analogue to \cite[Theorem~1.3]{KKL2022} and \cite[Theorem~1.2]{LiWangYi2023pre}. By \eqref{eqn_min_degree}, we can see that the smallest degree of ideals in any covering system of $O_S$ is bounded. By Theorem~\ref{mainthm}, it does not exist for sufficiently large $q$. Hence Theorem~\ref{mainthm} improves \eqref{eqn_min_degree} for sufficiently large $q$.
\end{remark}

For the proof of Theorem~\ref{mainthm}, we essentially follow the approaches in the works \cite{KKL2022, LiWangYi2023pre} and apply the distortion method. The main novelty in our proof is that we only need to use the estimate of the second moments in the distortion method. In \cite{KKL2022, LiWangYi2023pre}, the authors have to use the estimate of the first moments as well to control the summation from the part of small norms. 

In Section~\ref{sec_distortion_method}, we unify the distortion method in the settings of number fields and finite fields in \cite{LiWangYi2023pre, LiWangWangYi2023pre} to the setting of Dedekind domains. Then in Section~\ref{sec_pf_mainthm}, we apply an explicit upper bound for the number of prime ideals in $O_S$ to give an explicit estimate on the second moments. Then we apply the distortion method to give a proof of Theorem~\ref{mainthm}.

\section{The distortion method for Dedekind domains} \label{sec_distortion_method}
In this section, we illustrate the distortion method in the setting of Dedekind domains. The ring $O_S$ of $S$-integers considered in this paper and the ring of algebraic integers in a number field are two classical examples of Dedekind domains. Their nontrivial quotient rings are finite. We will also give an estimation of the second moments appearing in this method, which will be used in the proof of Theorem~\ref{mainthm}.

\subsection{The distortion method}
Let $R$ be a Dedekind domain. Assume that the quotient ring $R/I$ is finite for any nonzero ideal $I$ in $R$, and we define the norm of $I$ to be $|I|\colonequals |R/I|$. Then for any two nonzero ideals $I,J$ in $R$, we have $|IJ|= |I||J|$. Let $\mathcal{A}=\{A_{i}\colon 1\leq i\leq m\}$ be a finite collection of congruences of multiplicity $s$ in $R$, where $A_i=a_i+I_i, a_i\in R$ and $I_i$ is a nonzero ideal in $R$. Here the multiplicity of $\cA$ is defined to be $s=m(\cA)=\max_{ \text{Ideals } I \subset R}\#\{1\leq i\leq m\colon I_i=I\}$. Let $Q=I_1\cap\cdots\cap I_m$. Since $R$ is a Dedekind domain, we have a unique prime ideal decomposition of $Q$. We write it as $Q=\prod_{i=1}^J P_i^{\nu_i}$, where the $P_i$'s are distinct prime ideals and $|P_1|\leq |P_2|\leq \cdots \leq |P_J|$. Let $Q_j=\prod_{i=1}^jP_i^{\nu_i}, 1\leq j \leq J$ and $Q_0=R$.  For $1\le j \le J$, we put 
\[
\cB_j\colonequals\bigcup_{\substack{1\le i\le m\\ I_i|Q_j, I_i\nmid Q_{j-1}}} \set{a+Q\colon a \equiv a_i \mod{I_i}}.
\]

By the assumption on $R$, we have $R/Q_j$ is finite for each $j$. In particular, $R/Q=R/Q_J$ is finite. Let $\pi_j\colon R/Q \to  R/Q_j$ be the natural projection, and let $\pi_0$ be the trivial map, i.e., $\pi_0\colon R/Q \to \{0\}$. For any $a\in R$, we put $\bar{a}=a+Q\in R/Q$ and let 
\[
F_j(\bar{a})=\set{\bar{r}\in R/Q\colon \pi_j(\bar{r})=\pi_j(\bar{a})}
\]
be the fibre at $\bar{a}$. Observe that
\[
    \#F_j(\bar{a})=\#\ker \pi_j=\frac{|Q|}{|Q_j|}.
\]
For any $\bar{a}\in R/Q$, we define
\[
\alpha_j(\bar{a})=\frac{\# (F_{j-1}(\bar{a})\cap\cB_j)}{\#F_{j-1}(\bar{a})}.
\]

Now, we inductively define a sequence of probability measures $\P_0, \P_1, \dots, \P_J$ on $R/Q$. Let $\P_0$ be the uniform probability measure. Let $\delta_1,\dots,\delta_J\in[0,\frac12]$. Assuming we have defined $\P_{j-1}$, we define $\P_j$ as follows:
\[
\P_j(\bar{a})\colonequals \P_{j-1}(\bar{a}) \cdot \left\{
\begin{aligned}
	&\frac{1_{\bar{a}\notin \cB_j}}{1-\alpha_j(\bar{a})} & \text{if }  \alpha_j(\bar{a})<\delta_j,\\
	&\frac{\alpha_j(\bar{a})-1_{\bar{a}\in \cB_j}\delta_j}{\alpha_j(\bar{a})(1-\delta_j)} & \text{if }  \alpha_j(\bar{a})\geq \delta_j.
\end{aligned}
\right.
\]

For $k\in\Z_{\geq 1}$, we define the $k$-th moment of $\alpha_j$ as follows:
\begin{equation}\label{eqn_kth_moment_definition}
    M_j^{(k)}\colonequals\sum_{\bar{a}\in R/Q}\alpha_j^k(\bar{a})\P_{j-1}(\bar{a}).
\end{equation}

By \cite[Theorem~2.1]{LiWangYi2023pre}, if 
	\[
	\sum_{j=1}^J
	\min\left\{  M_j^{(1)},\frac{M_j^{(2)}}{4\delta_j(1-\delta_j)}\right\}<1 ,
	\]
	then 
 \[
 \P_J\of{\bigcup_{1\le j \le J} \cB_j}<1
 \]
 and hence $\cA$ cannot cover $R$. Thus, the following theorem holds if we take $\delta_j=\frac{1}{2}$ for all $1\leq j \leq J$. It is a variant of  the distortion method introduced in  Balister et al.'s work \cite[Theorem~3.1]{BBMST2022} and Klein et al.'s work \cite{KKL2022}.
\begin{theorem}\label{thm_distortion_method}
With the notation as above,	take $\delta_j=\frac{1}{2}$ for all $1\leq j \leq J$. If
	\[
	\sum_{j=1}^J M_j^{(2)}<1,
	\]
	then $\cA$ does not cover $R$. 
\end{theorem}

\subsection{Bounding the second moments} In this subsection, we give an explicit upper bound for $M_j^{(2)}$ in Lemma~\ref{lem_second_moments_estimate}. This is extracted from the works \cite{KKL2022, LiWangYi2023pre}. We provide a proof here for reader's convenience. First, we recite two results that are generalizations  of  \cite[Lemma~3.2]{KKL2022}, \cite[Lemma~3.4]{BBMST2022}, and \cite[Lemmas~3.1 and 3.2]{LiWangYi2023pre}. For the proofs of the following Lemmas~\ref{lem_alpha_estimate} and \ref{lem_p_j_estimate}, one may follow the proofs in \cite{BBMST2022,KKL2022,LiWangYi2023pre}.

\begin{lemma}\label{lem_alpha_estimate}
	For any $\bar{a}\in R/Q$ and $1\leq j \leq J$, we have
\[
\alpha_j(\bar{a})\leq \sum_{r=1}^{\nu_j} \sum_{H\mid Q_{j-1}} \sum_{\substack{1\leq i \leq m\\ I_i=HP_j^r}} \frac{1_{\bar{a}\subseteq a_i+H}}{|P_j|^r}.
\]
\end{lemma}

\begin{lemma}\label{lem_p_j_estimate}
For each $0\leq j\leq J, a\in R$, and for any ideal $I\subseteq R$ such that $I \mid Q$, we have
\[
\P_j(a+I)\leq \frac{1}{|I|}\prod_{P_i\mid I,\ i\leq j}\frac{1}{1-\delta_i}.
\]
\end{lemma}

Now, following the proofs of \cite[Lemma~3.3]{KKL2022} and \cite[Lemma~3.3]{LiWangYi2023pre}, we will use Lemmas~\ref{lem_alpha_estimate} and \ref{lem_p_j_estimate} to estimate $M_j^{(2)}$ explicitly.

\begin{lemma}\label{lem_second_moments_estimate}
Let $s=m(\cA)$ and $\delta_1,\dots,\delta_J\in[0,\frac12]$. We have
    \[
     M_j^{(2)} \leq \frac{s^2}{(|P_j|-1)^2}\prod_{i<j}\left(1+\sum_{\nu\geq 1}\frac{4\nu+2}{|P_i|^\nu} \right)
    \]
for all $1\leq j\leq J$.
\end{lemma}

\begin{proof} By definition, $M_j^{(2)}=\sum_{\bar{a}\in R/Q}\alpha_j^2(\bar{a})\P_{j-1}(\bar{a})$. So by Lemma~\ref{lem_alpha_estimate} we have
\begin{equation}\label{eqn_2nd_moment}
    M_j^{(2)}\leq \sum_{1\leq r_1, r_2\leq \nu_j} \sum_{H_1,H_2 \mid Q_{j-1}} \sum_{\substack{1\leq i_1, i_2\leq m\\ I_{i_\ell}=H_\ell P_j^{r_\ell},\,\ell=1,2}} \frac{\P_{j-1}\left(\bigcap_{\ell=1}^2(a_{i_\ell}+H_\ell)\right)}{|P_j|^{r_1+r_2}}.
\end{equation}
We observe that given $r_1, r_2$ and $H_1, H_2$, there are at most $s^2$ choices for $i_1, i_2$ such that $I_{i_\ell}=H_\ell P_j^{r_\ell}$. By the Chinese remainder theorem, for each such choice of $i_1, i_2$, the set $\bigcap_{\ell=1}^2(a_{i_\ell}+H_\ell)$ is either empty or a congruence class modulo $H_1\cap H_2$. Moreover, by Lemma~\ref{lem_p_j_estimate} we get that 
\begin{equation}\label{eqn_prob_of_intersection_of_APs}
    \P_{j-1}\left(\bigcap_{\ell=1}^2(a_{i_\ell}+H_\ell)\right) \leq \frac{1}{|H_1\cap H_2|}\prod_{\substack{i\leq j-1\\ P_i\mid H_1\cap H_2}}(1-\delta_i)^{-1}
\end{equation}
for at most $s^2$ choices of $1\leq i_1, i_2\leq m$.
It follows by plugging \eqref{eqn_prob_of_intersection_of_APs} into \eqref{eqn_2nd_moment} that
\[
M_j^{(2)}\leq s^2 \sum_{1\leq r_1, r_2\leq \nu_j} \sum_{H_1, H_2 \mid Q_{j-1} } \frac{1}{|H_1\cap H_2|\cdot|P_j|^{r_1+r_2}}\cdot\prod_{\substack{i\leq j-1\\ P_i\mid H_1 \cap H_2 }}(1-\delta_i)^{-1}.
\]

Since $\delta_i\in[0,1/2]$ for all $1 \leq i \leq J$, we have
\begin{equation*}
    \prod_{P_i\mid H_1 \cap H_2}(1-\delta_i)^{-1}\leq 2^{\omega(H_1 \cap H_2)}.
\end{equation*}
Here, $\omega(H)$ denotes the number of distinct prime ideal factors of $H$ for any ideal $H\subseteq R$. Notice that $2^{\omega(H)}$ is a multiplicative function for $H\subseteq R$. Hence we obtain
\begin{align*}
   M_j^{(2)} & \leq s^2\sum_{1\leq r_1, r_2\leq \nu_j}\sum_{H_1, H_2\mid Q_{j-1}}\frac{2^{\omega(H_1 \cap H_2)}}{|H_1 \cap H_2| \cdot |P_j|^{r_1+r_2}}\\
   &\leq \frac{s^2}{(|P_j|-1)^2}\sum_{H_1, H_2\mid Q_{j-1}}\frac{2^{\omega(H_1 \cap H_2)}}{|H_1 \cap H_2|}\\
   &\leq \frac{s^2}{(|P_j|-1)^2}\prod_{i<j}\left(1+\sum_{\nu\geq 1}\frac{2}{|P_i|^\nu} (2\nu+1)\right),
\end{align*}
as desired.  
\end{proof}

\section{Proof of Theorem~\ref{mainthm}}\label{sec_pf_mainthm}

In this section we adopt the notation of Section~\ref{sec_distortion_method} and take $R=O_S$. We will prove Theorem~\ref{mainthm} by applying Theorem~\ref{thm_distortion_method} and Lemma~\ref{lem_second_moments_estimate}. The numerical estimates on the infinite series in the proofs are done by the Mathematica software \cite{Mathematica}. Let $\pi_{K, S}(n)$ be the number of prime ideals of degree $n$ in $O_S$. By the analogue of the prime number theorem over global function fields (see \eqref{eqn_PNT_for_FF}), we have that
$$\pi_K(n)=\frac{q^n}{n}+O\left(\frac{q^{n/2}}{n}\right).$$
Moreover, by the proof of \cite[Theorem~5.12]{Rosen2002}, we can see that
\begin{equation}\label{eqn_PPT_explicit_upper_bound}
    \pi_{K, S}(n)\leq \pi_K(n)\leq \frac{q^n}{n}+\frac{(2+4g)q^{n/2}}{n}
\end{equation}
for $q\geq4$, where $g$ is the genus of $K$.  To estimate $M_j^{(2)}$, we first use \eqref{eqn_PPT_explicit_upper_bound} to estimate the sum of reciprocals of the norms of prime ideals in $O_S$ as follows. 

\begin{lemma}\label{lem_mertens_estimates}
For $q\geq 700$, we have
\[
	\sum_{\deg P\leq N}\frac{1}{|P|}\leq \log N+1.08+0.16g,  \quad  \sum_{P}\frac{1}{|P|^2}\leq \frac{1.08+0.16g}{q},
\]
where the term $P$ runs through all prime ideals in $O_S$.
\end{lemma}
\begin{proof} For the first estimate, by \eqref{eqn_PPT_explicit_upper_bound} and the following inequality 
$$\sum_{n=1}^N \frac{1}{n}\leq \log N+1,$$
we have that
\begin{align*}
    \sum_{\deg P\leq N}\frac{1}{|P|}&=\sum_{n=1}^N\frac{\pi_{K, S}(n)}{q^n}\\
    &\leq \sum_{n=1}^N \frac{1}{n} + \sum_{n=1}^N\frac{2+4g}{nq^{n/2}}\\
    &<\log N+1+ \sum_{n=1}^\infty\frac{2+4g}{n\cdot 700^{n/2}}\\
    &<\log N +1.08+0.16g.
\end{align*}
Here, $\sum_{n=1}^\infty\frac{2}{n\cdot 700^{n/2}}=0.0770585\cdots$.

For the second estimate, we have
\begin{align*}
    \sum_{P}\frac{1}{|P|^2}&=\sum_{n=1}^\infty\frac{\pi_{K, S}(n)}{q^{2n}}\\
    & \leq \sum_{n=1}^\infty\frac{1}{q^{2n}}\left(\frac{q^n}{n}+\frac{(2+4g)q^{n/2}}{n}\right)\\
    & =\frac1q \sum_{n=1}^\infty\left(\frac{1}{nq^{n-1}}+\frac{2+4g}{nq^{3n/2-1}}\right)\\
    & \leq \frac{1}{q} \sum_{n=1}^\infty \left(\frac{1}{n\cdot 700^{n-1}}+\frac{2+4g}{n\cdot 700^{3n/2-1}}\right)\\
    &< \frac{1.08+0.16g}{q}.
\end{align*}
Here, $\sum_{n=1}^\infty \left(\frac{1}{n\cdot 700^{n-1}}+\frac{2}{n\cdot 700^{3n/2-1}}\right)=1.07631\cdots$ and $\sum_{n=1}^\infty \frac{4}{n\cdot 700^{3n/2-1}}=0.15119\cdots$.
\end{proof}

Now, we derive upper bounds for $M_j^{(2)}$ explicitly, where $1 \leq j \leq J$, by utilizing Lemma~\ref{lem_second_moments_estimate} and Lemma~\ref{lem_mertens_estimates}.

\begin{lemma}\label{lem_2nd_moments_estimate}
Suppose $q\geq 700$. Let $s=m(\cA)$ and $\delta_1,\dots,\delta_J\in[0,\frac12]$. Then we have that
    \[
     M_j^{(2)} \leq \frac{1.01e^{6.5+0.97g}s^2(\deg P_j)^6}{|P_j|^2}
    \]
for all $1\leq j\leq J$.
\end{lemma}

\begin{proof} By Lemma~\ref{lem_second_moments_estimate}, we have
\begin{align*}
   M_j^{(2)} &\leq \frac{s^2}{(|P_j|-1)^2}\prod_{i<j}\left(1+\sum_{\nu\geq 1}\frac{2}{|P_i|^\nu} (2\nu+1)\right)\\
   &=\frac{s^2}{(|P_j|-1)^2}\prod_{i<j}\left(1+\frac{6}{|P_i|}+\frac{10|P_i|-6}{|P_i|(|P_i|-1)^2}\right)\\
   &\leq \frac{1.01s^2}{|P_j|^2}\prod_{i<j}\left(1+\frac{6}{|P_i|}+\frac{11}{|P_i|^2}\right)\\
   &\leq \frac{1.01s^2}{|P_j|^2}\,\mathrm{exp}\left\lbrace \sum_{i<j}\frac{6}{|P_i|}+\sum_{i<j}\frac{11}{|P_i|^2}\right\rbrace
\end{align*}
by the inequality $1+t\leq e^t$ for all $t\in\R$. Then by Lemma~\ref{lem_mertens_estimates}, we get that 
\begin{align*}
   M_j^{(2)}&\leq  \frac{1.01s^2}{|P_j|^2}\,\mathrm{exp} \left\lbrace 6(\log \deg P_j+ 1.08+0.16g) + 11\cdot (1.08+0.16g)/700 \right\rbrace\\
   &<  \frac{1.01s^2}{|P_j|^2}\,\mathrm{exp} \left\lbrace 6\log \deg P_j +6.5+0.97g \right\rbrace\\
   &= \frac{1.01e^{6.5+0.97g}s^2(\deg P_j)^6}{|P_j|^2},
\end{align*}
as desired.
\end{proof}

\noindent\textit{Proof of Theorem~\ref{mainthm} and Corollary~\ref{maincor}.} With the notation as above, using the estimate of the second moment $M_j^{(2)}$ as in Lemma~\ref{lem_2nd_moments_estimate}, we have
\begin{align*}
	\sum_{1\leq j\leq J}M_j^{(2)} & \leq \sum_{1\leq j\leq J} \frac{1.01e^{6.5+0.97g}s^2(\deg P_j)^6}{|P_j|^2}\\
	&<  \sum_{P} \frac{1.01e^{6.5+0.97g}s^2(\deg P)^6}{|P|^2}\\
	&=  1.01e^{6.5+0.97g}s^2 \cdot \sum_{n=1}^\infty \frac{n^6}{q^{2n}} \pi_{K, S}(n)\\
	&\leq 1.01e^{6.5+0.97g}s^2 \cdot \sum_{n=1}^\infty \frac{n^6}{q^{2n}}\left(\frac{q^n}{n}+\frac{(2+4g)q^{n/2}}{n}\right)\\
	&\leq  \frac{1.01e^{6.5+0.97g}s^2}{q} \sum_{n=1}^\infty \left(\frac{n^5}{700^{n-1}} + \frac{(2+4g)n^5}{700^{3n/2-1}}\right)\\
	&< \frac{(1.14+0.16g)e^{6.5+0.97g}s^2 }{q}
\end{align*}
for $q\ge 700$. Here, $\sum_{n=1}^\infty \left(\frac{n^5}{700^{n-1}} + \frac{2n^5}{700^{3n/2-1}}\right)=1.12194\cdots$ and $\sum_{n=1}^\infty \frac{4n^5}{n\cdot 700^{3n/2-1}}=0.151447\cdots$. 

Thus, if 
\[
q\geq (1.14+0.16g)e^{6.5+0.97g}s^2,
\]
then we have 
\[
\sum_{1\leq j\leq J}M_j^{(2)}<1.
\]
By Theorem~\ref{thm_distortion_method}, we conclude that $\cA$ is not a covering system for $O_S$ for $q\geq (1.14+0.16g)e^{6.5+0.97g}s^2$ and Theorem~\ref{mainthm} follows. Then Corollary~\ref{maincor} follows immediately by taking $g=0$ and $1.14e^{6.5}=758.261\cdots$ in Theorem~\ref{mainthm}. This completes the proof.\qed

\section*{Acknowledgments}
The authors would like to thank the referee for a very careful review and a number of precise comments and helpful suggestions, which improve this paper a lot. Huixi Li's research is partially supported by the National Natural Science Foundation of China (Grant No.~12201313). Shaoyun Yi is supported by the National Natural Science Foundation of China (No.~12301016) and the Fundamental Research Funds for the Central Universities (No.~20720230025).

\end{document}